\documentclass[11pt]{article}
\usepackage{epsfig}
\usepackage{verbatim}
\usepackage{amsmath}
\usepackage[pdfmenubar=false, letterpaper=false, pdfpagelabels=true]{hyperref}
\usepackage{makeidx,showidx}
\usepackage{pst-node}  

\begin{document}

\newcommand{\bm}[1]{\mbox{\boldmath$#1$}}
\def\mvec#1{{\bm{#1}}}   

\title{On the Peirce's {\it balancing reasons rule} failure\\
in his ``large bag of beans'' example}

\author{G.~D'Agostini \\
Universit\`a ``La Sapienza'' and INFN, Roma, Italia \\
{\small (giulio.dagostini@roma1.infn.it,
 \url{http://www.roma1.infn.it/~dagos})}
}

\date{}

\maketitle

\begin{abstract}
Take a large bag 
of black and white beans, with all possible 
proportions considered initially equally likely,
and imagine to make random extractions with reintroduction.
Twenty consecutive observations of black make us
highly confident that the next bean will be black too.
On the contrary, the observation of 
1010 black beans and 990 white ones leads us to 
judge the two possible outcomes about equally probable. 
According to C.S. Peirce this reasoning 
violates what he called ``rule of {\it balancing reasons}'',
because the difference of ``arguments'' in favor and against
the outcome of black is 20 in both cases. Why? (I.e. why does 
that rule not apply here?)    
\end{abstract}

\section{Introduction}
Let us take the following example from C.S. Peirce's
{\it The probability of induction}\cite{Peirce}:
\begin{quote}
{\sl 
``Suppose we have a large bag of beans from which one has
been secretly taken at random and hidden under a thimble.
We are now to form a probable judgement of the color of that
bean, by drawing others singly from the bag and looking at them, 
each one to be thrown back, 
and the whole well mixed up after each drawing. \\
\ldots \\
Suppose that the first bean which we drew from our bag
were black. That would constitute and argument, no matter
how slender, that the bean under the thimble was also black.
If the second bean were also to turn out black, that would be
a second independent argument re\"enforcing the first.
If the whole of the first twenty beans drawn should prove black, 
our confidence that the hidden bean was black would justly attain 
considerable strength. But suppose the twenty-fits bean were to
be white and that we were to go on drawing until we found that we had
drawn 1,010 black beans and 990 white ones. 
We would conclude that our first twenty beans being black was 
simply an extraordinary accident, and that in fact
the proportion
of white beans to black was sensible equal, and
that it was an even chance that the hidden bean was black.
Yet according to the rule of {\it balancing reasons},
since all the drawings of black beans are so many independent
arguments in favor of the one under the thimble being black,
and all the white drawings so many against it, 
an excess of twenty black beans ought to produce the same degree
of belief that the hidden bean was black, 
whatever the total number drawn.''
}\cite{Peirce}
\end{quote}
The philosopher does not try to resolve the manifest contradiction
in the rest of the article and the question is then left
to the reader as a kind of paradox of what he calls
the ``conceptualistic view of probability'' (nowadays `subjective
probability'), although its solution is rather easy: 
his `rule of {\it balancing reasons}' does not apply 
to the first practical example he provides, because 
the `arguments' are not independent. 

\section{Which box? Which color?}
Let us think to a slight different problem. 
We have two boxes, $B_1$ and $B_2$, 
containing well known proportions $p_1$ and $p_2$ 
of white balls, respectively (the remaining one are black). 
If we make random extractions ($E_i$)
with reintroduction, the probability of getting black (B)
and white (W) balls are: 
\begin{eqnarray}
P(E_i=\mbox{W}\,|\,B_1,I) &=& p_1 \\
P(E_i=\mbox{B}\,|\,B_1,I) &=& 1-p_1 \\
P(E_i=\mbox{W}\,|\,B_2,I) &=& p_2 \\
P(E_i=\mbox{B}\,|\,B_2,I) &=& 1-p_2\,, 
\end{eqnarray}
where the symbol `$|$' stands for `given', i.e. 'under the condition',
whereas the  ubiquitous `$I$' stands for the
state of information under which probability values are assessed.

If we take one of the boxes at random (hereafter $B_?$) 
this could be equally likely $B_1$ or $B_2$ and then 
the probability of getting black or white will be the averages
of the probabilities given the two box compositions. 
As soon as we start sampling the box content by extractions
followed by reintroduction our opinion 
concerning the box composition is modified
by the experimental information, and the 
probability of occurrence of 
white in the next extraction is modified too.\footnote{If 
the amount of balls in the box is very large, 
thinking to the next extraction, or taking at random 
a ball at the very beginning and hiding it ``under a thimble'',
is practically the same.}
\subsection{Probability of white/black from a box on unknown composition}
In general, if, after $n$ observations, our beliefs in 
the two kinds of boxes are $P(B_1\,|\,\mvec{E},B_?,I)$
and $P(B_2\,|\,\mvec{E},B_?,I)$, 
the probability that a next extraction gives white is given by
\begin{eqnarray}
P(\mbox{W}\,|\,\mvec{E},B_?,I) 
&=& P(\mbox{W}\,|\,B_1,I)\,P(B_1\,|\,\mvec{E},B_?,I) + 
     P(\mbox{W}\,|\,B_2,I)\,P(B_2\,|\,\mvec{E},B_?,I) \label{eq:W_E_PBi}\\
&=& p_1\,P(p_1\,|\,\mvec{E},B_?,I) +
    p_2\,P(p_2\,|\,\mvec{E},B_?,I)\,, \label{eq:W_E_Ppi}
\end{eqnarray}
where $\mvec{E}$ is the ensemble of all observations, 
i.e. $\mvec{E}=\{E_1,E_2,\ldots,E_n\}$. In Eq.~(\ref{eq:W_E_Ppi})
$P(B_i\,|\,\mvec{E},B_?,I)$ has been replaced by 
$P(p_i\,|\,\mvec{E},B_?,I)$ to remark that our belief 
in a given proportion is equal to the our belief in the
corresponding box type. Note that, since 
$P(p_1\,|\,\mvec{E},B_?,I)+P(p_2\,|\,\mvec{E},B_?,I)=1$, 
Eq.~(\ref{eq:W_E_Ppi}) can be read as weighted average. 
If, instead of just two possible box compositions, we have
many, Eq.~(\ref{eq:W_E_Ppi}) becomes
\begin{eqnarray}
P(\mbox{W}\,|\,\mvec{E},B_?,I) 
&=& \sum_i p_i\,P(p_i\,|\,\mvec{E},B_?,I)\,. \label{eq:W_E_Ppi_many}
\end{eqnarray}
\subsection{Probability of the different box composition given the past observations}
As far as the updating of probability is concerned,  
the most convenient way in the case
of two hypotheses is to use the 
update of probability ratios (odds) via the Bayes factor. 
Since the events  black or white are independent 
{\it given a box compositions} 
and using the notation of 
Ref.~\cite{Colombo}\footnote{This paper is strictly related
to Ref.~\cite{Colombo}, because I discovered Peirce's 
{\it The Probability of Induction} making a short historical 
research on the use of the logarithmic updating of
odds (see Appendix E there).} 
(sections 2.3 and
2.4), we can write
\begin{eqnarray}
O_{1,2}(\mvec{E},I) &=& \tilde O_{1,2}(\mvec{E},I) \times O_{1,2}(I)\,, 
\label{eq:BF}
\end{eqnarray}
where  the priors odds $ O_{1,2}(I)$
are unitary in this case ('even'), while the overall
Bayes factor is
\begin{eqnarray}
\tilde O_{1,2}(\mvec{E},I) &=& \prod_{k=1}^n
\tilde O_{1,2}(E_k,I)\,,  \label{eq:product_Odds_1}
\end{eqnarray}
with 
\begin{eqnarray}
O_{1,2}(\mvec{E},I) &=& \frac{P(p_1\,|\,\mvec{E},B_?,I)}
                             {P(p_2\,|\,\mvec{E},B_?,I)} \\
\tilde O_{1,2}(E_k,I) &=& \frac{P(E_k\,|\,p_1,I)}{P(E_k\,|\,p_2,I)}\,,
\end{eqnarray}
where the Bayes factors due to each piece of evidence are written
as $\tilde O_{1,2}(E_k,I)$ to remark that they {\it would be}
the odds only considering the individual piece of evidence 
$E_k$, provided the two hypotheses were otherwise considered
equally likely. 

\subsubsection{Logarithmic update and weight of evidence}
The update rule (\ref{eq:BF}) can be turned into 
an additive rule if,  as first (as far as I know) 
proposed by Peirce
in the same paper of the bag of beans example, we 
take the logarithm of it. Using the notation 
of Ref.~\cite{Colombo}, 
we can rewrite Eq.~(\ref{eq:BF}) 
as
\begin{eqnarray}
\mbox{JL}_{1,2}(\mvec{E},I) &=& 
 \Delta\mbox{JL}_{1,2}(\mvec{E},I) + \mbox{JL}_{1,2}(I)
\,, 
\label{eq:JLupdate}
\end{eqnarray}
with 
\begin{eqnarray}
 \Delta\mbox{JL}_{1,2}(\mvec{E},I)  &=& 
\sum_{k=1}^n \Delta\mbox{JL}_{1,2}(E_k,I)\,,
  \label{eq:Sum_DeltaJL_k}
\end{eqnarray}
where the JL's and the $\Delta$JL's are base 10 logarithms of
odds and of Bayes factors, respectively. 
The JL's ({\it judgement leaning}) correspond to Peirce's
{\it intensities of belief}, their variation ($\Delta$JL) being
due to the {\it weight of evidence} (see Ref.~\cite{Peirce} and 
Appendix E of Ref.~\cite{Colombo}). 

The contributions $\Delta\mbox{JL}_{1,2}(E_k,I)$ can be 
positive or negative, depending if the corresponding
Bayes factors are larger or smaller than one, 
and they are considered by Peirce as {\it arguments} in favor
or against the hypothesis 1 ($B_1$, or $p_1$, here):
\begin{quote}
{\sl
``The rule of the combination of independent arguments takes
a very simple form when expressed in terms of the intensity
of belief, measured in the proposed way. It is this:
Take the sum of all the feelings of belief which would be produced
separately by all the arguments {\it pro}, subtract from that the 
similar sum for arguments {\it con}, and the remainder is the feeling
of belief which we ought to have on the whole.
This is a proceeding which 
men often resort to, 
under the name of {\it balancing reasons}.''
}\cite{Peirce}
\end{quote}
At this point we only need to write down the weights of evidence
due to the observation of the different colors:
\begin{eqnarray}
\Delta\mbox{JL}_{1,2}(E_k=\mbox{W},I) &=& \log_{10}\frac{p_1}{p_2} \\
\Delta\mbox{JL}_{1,2}(E_k=\mbox{B},I) &=& \log_{10}\frac{1-p_1}{1-p_2}\,.
\end{eqnarray}
As we we see, the absolute weight of `arguments' depends 
on the values of $p_1$ and $p_2$. If they are very similar,
the indication provided by the experimental information 
is very week and we need a very large number of observations
to discriminate between the two hypotheses
(see e.g. Appendix G of Ref.~\cite{Colombo}). If, instead,
the proportion of one kind of balls is very close to 
zero or to 1, the indications can be rather strong
and just one or a few extractions make us highly confident
about the box composition. At the limit, 
if one of the box only contains white or black balls, a single
observation showing the opposite 
colors is enough to rule out that hypothesis 
($|\Delta\mbox{JL}_{1,2}(E_k,I)|=\infty$).

\subsubsection{Combined weight of evidence and final odds
after a sequence of extractions}
Since the weights of evidence due to independent pieces of 
evidence sum up, after the observation of $n_W$ white and $n_B$ black
we have 
\begin{eqnarray}
\Delta\mbox{JL}_{1,2}(n_W,n_B,I) &=& n_W\log_{10}\frac{p_1}{p_2} 
+ n_B \log_{10}\frac{1-p_1}{1-p_2}\,,
\end{eqnarray}
or, since we started from uniform priors,
\begin{eqnarray}
O_{1,2}(n_W,n_B,I) &=& \left(\frac{p_1}{p_2}\right)^{n_W}\,\left(\frac{1-p_1}{1-p_2}\right)^{n_B} \\
 &=& \frac{p_1^{n_W}\,(1-p_1)^{n_B}}
          {p_2^{n_W}\,(1-p_2)^{n_B}}
\,.
\end{eqnarray}
The probabilities of the two box compositions are then 
\begin{eqnarray}
P(p_i\,|\,n_W,n_B,I)  &=&  \frac{p_i^{n_W}\,(1-p_i)^{n_B}}
   {p_1^{n_W}\,(1-p_1)^{n_B} + p_2^{n_W}\,(1-p_2)^{n_B}}\,.
\end{eqnarray}
This formula can easily extended to the case of many box composition:
\begin{eqnarray}
P(p_i\,|\,n_W,n_B,I)  &=&  \frac{p_i^{n_W}\,(1-p_i)^{n_B}}
   {\sum_j p_j^{n_W}\,(1-p_j)^{n_B}}\,. \label{eq:many_p}
\end{eqnarray}

\subsection{Case with two symmetric bag compositions}
A particular case, that can be useful to clarify 
the difference with respect to the different problem 
 discussed in the following section,
is when $p_2=1-p_1$, for example $p_1=1/4$ and $p_2=3/4$. 
In this case we have 
\begin{eqnarray}
\Delta\mbox{JL}_{1,2}(E_k=\mbox{W},I) &=&  \log_{10}\frac{p_1}{1-p_1} \\
\Delta\mbox{JL}_{1,2}(E_k=\mbox{B},I) &=&  \log_{10}\frac{1-p_1}{p_1}
\end{eqnarray}
i.e. 
\begin{eqnarray}
\Delta\mbox{JL}_{1,2}(E_k=\mbox{B},I) &=& 
 -\Delta\mbox{JL}_{1,2}(E_k=\mbox{W},I)
\end{eqnarray}
the weights of evidence provided by black and white have opposite
sign but are equally in module. It follows that our judgement in favor
of the two boxes depends only on the difference of black and white
balls observed, but not on the number of extractions. Here then 
the rule of balancing reasons applies.

\section{From many to (virtually) infinite box compositions}
\begin{figure}
\centering\epsfig{file=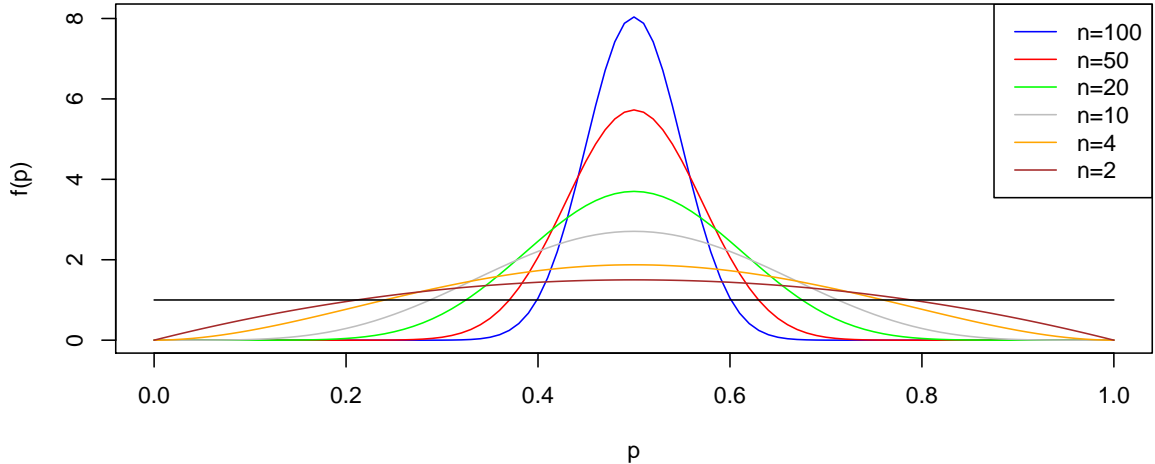,clip=,width=\linewidth}
\caption{\sl Probability density function of the white ball proportion
$p$, having observed 
exactly 50\% white balls in 2, 4, 10, 20, 50 
and 100 extractions (the curves become in order narrower).
The horizontal line represents
the uniform prior.}
\label{fig:beta_n}
\end{figure}

In the limit that the number of possible compositions 
is virtually infinite, the parameter $p$ that gives the white
ball proportion becomes continuous and the problem is 
solved in terms of  probability 
density function  $f(p\,|\,\mvec{E},I)$. Essentially
Eq.~(\ref{eq:many_p}) becomes 
\begin{eqnarray}
f(p\,|\,n_W,n_B,I)  &=&  \frac{p^{n_W}\,(1-p)^{n_B}}
   {\int_0^1 p^{n_W}\,(1-p)^{n_B}\,dp}\,, \label{eq:infinite_p}
\end{eqnarray}
that gives
\begin{eqnarray}
f(p\,|\,n_W,n_B,I)  &=&  \frac{(n+1)!}{n_W!\,n_B!}\,\,p^{n_W}
                         (1-p)^{n_B}\,. \label{eq:f(p)}
\end{eqnarray}
Some examples of $f(p\,|\,n_W,n_B,I)$ are given in figure 
\ref{fig:beta_n} with for several numbers of extractions, 
assuming that in all cases the fraction of white balls
has been 50\%.
We see that with the increasing number of extractions
we get more and more confident that $p$ is around 
0.5.\footnote{The uncertainty, measured by the standard 
deviation of the distribution $\sigma(p)$, is given 
by
\begin{eqnarray*}
\sigma(p\,|\,n_W,n_B,I) &=& \sqrt{E[p](1-E[p])\,\frac{1}{n_W+n_B+3}}
\propto \frac{1}{\sqrt{n+3}}\,, \label{eq:pdf_p}
\end{eqnarray*}
with $E[p]$ equal to the {\it expected value}, given 
by Eq.~(\ref{eq:Int_p}).\\
Note that $\sigma(p)$ is the uncertainty about the proportion of 
white balls in the box and not about the probability of
having white in the next extraction, which is {\it exactly} 1/2
in all cases in which an equal number of black and white balls
has been observed! 
It seems that this point was not very clear to Peirce,
who writes in Ref.~\cite{Peirce}, also referring to the 
large ``bag of bean'', after the first period
of the quote reported in page 1: 
\begin{quote}
{\sl ``Suppose the first drawing is white and the next black.
We conclude that there is not an immense preponderance of either 
color, and that there is something like an even chance that the bean
under the thimble is black. But this judgement might be altered
by the next few drawings. When we
have drawn ten times, if 4, 5, or 6, are white,
we have more confidence that the chance [{\rm  that the bean
under the thimble is black, we have to understand}] is even.
When we have drown a thousand
times, if about half have been white, we have great 
confidence in this result.''
}\cite{Peirce}
\end{quote}
To be more precise, there are several things that should kept
separate in our reasonings:
\begin{itemize}
\item The proportion of white balls in the box, that is $p$, our
      uncertainty concerning it being described by the probability
      density function (\ref{eq:f(p)}), with expected value 
      given by Eq.~(\ref{eq:Int_p}) and `standard uncertainty' 
      $\sigma(p)$ given at the beginning of this footnote 
      (results valid from a uniform prior).
\item The relative frequency of white balls 
      that we expect in a series of $m$ extraction, given the 
      past observation of $n_W$ and $n_B$. The expression that gives
      our beliefs on all possible (in number $m+1$) values of 
      the relative
      frequency is quite complicate and can be found in 
      section 7.3 of Ref.~\cite{BR}. The case in which $m$ tends to infinite
      is instead rather easy to understand, since, calling 
      $\varphi_m$ the possible value 
      of the relative frequency in $m$
      extractions, we have, under the assumption that $p$ 
      is perfectly known,
      \begin{eqnarray*}
      \mbox{E}[\varphi_m\,|\,p,I] &=& p \\
      \sigma[\varphi_m\,|\,p,I] &=& \frac{\sqrt{p\,(1-p)}}{\sqrt{m}}\,.
      \end{eqnarray*}
      That is, in the limit $m\rightarrow\infty$, we feel practically
      sure to observe a value of $\varphi_\infty$ equal to $p$
      (`Bernoulli theorem'). 
      If we are, instead, uncertain about $p$, then we are uncertain about
      $\varphi_\infty$ exactly in the same way and the probability
      density function $f(\varphi_\infty\,|\,n_W,n_B,I)$ has the
      same shape of  $f(p\,|\,n_W,n_B,I)$. 
\item The probability that the next outcome will be white, the evaluation 
      of which has to take in consideration all possible values 
      of $p$, each weighted by how much we believe it (to be precise,
      since the proportion is virtually a real number, 
      the beliefs concern small intervals of $p$). But this is exactly 
      the expected value of $p$, according to the relation 
      (\ref{eq:E[p]}).
      In the particular case of {\it absolute symmetry} of our observations
      {\it and} of our prior, there should be {\it not the slightest
      rational preference} in favor of either color 
      [$P(\mbox{W}\,|\,n_W=n_B,B_?,I)=1/2$], no matter if
      we are very uncertain about box composition or future 
      relative frequencies. 
\end{itemize}
} 
 
The most probable values of $p$ are those around $n_W/(n_W+n_B)$,
although the probability to have a white ball in a future
observation has to take into account all possible
compositions in a manner similar to that seen in 
Eq.~(\ref{eq:W_E_Ppi_many}), i.e. 
\begin{eqnarray}
P(\mbox{W}\,|\,\mvec{E},B_?,I) 
&=& \int_0^1\! p\,f(p\,|\,\mvec{E},B_?,I)\,dp \,. \label{eq:W_E_fp}
\end{eqnarray}
In the r.h.s. of (\ref{eq:W_E_fp}) we recognize the {\it expected value}
of $p$, `barycenter' of the probability density function. 
Therefore it follows that 
\begin{eqnarray}
P(\mbox{W}\,|\,\mvec{E},B_?,I) &=& 
\mbox{E}[p\,|\,\mvec{E},B_?,I] = 
 \int_0^1\! p\,f(p\,|\,\mvec{E},B_?,I)\,dp \,. \label{eq:E[p]}
\end{eqnarray} 
The result of the integral (\ref{eq:W_E_fp})
is 
\begin{eqnarray}
\int_0^1\! p\,f(p\,|\,\mvec{E},B_?,I)\,dp &=& \frac{n_W+1}{n_W+n_B+2}\,,
\label{eq:Int_p}
\end{eqnarray} 
thus leading to 
the famous (although often misused!) {\it Laplace
rule of successions}
\begin{eqnarray}
P(\mbox{W}\,|\,\mvec{E},B_?,I) = \frac{n_W+1}{n_W+n_B+2}\label{P_W_E_P}
\end{eqnarray}
and, by symmetry,
\begin{eqnarray}
P(\mbox{B}\,|\,\mvec{E},B_?,I) = \frac{n_B+1}{n_W+n_B+2}\,.
\end{eqnarray}
Using Peirce numbers, we get
\begin{eqnarray}
P(\mbox{B}\,|\,20\mbox{B},B_?,I) &=& \frac{21}{22} = 95.5\% \\
P(\mbox{B}\,|\,1010\mbox{B}\!+\!990\mbox{W},B_?,I) &=& \frac{991}{2002} 
= 50.4\%\,, \label{eq:PB_1010_990}
\end{eqnarray}
providing quite different degrees of belief, as it is intuitive 
and as it was clear to Peirce, who, by the way, uses
Laplace rule of successions several times in Ref.~\cite{Peirce}.

\section{Weights of evidence in favor of a black or white 
bean under the thimble}
We have now all the tools to analyze Peirce's bag of bean example,  
in which the `arguments' did not regard the box composition,
but the occurrence of a white or black bean in a future
extraction (the fact that the bean
was extracted at the beginning is irrelevant, as it has 
already been observed). 

Since  the bag is `large' the proportion of white beans
can be considered as a real number $p$ ranging between 0 and 1.
Moreover, as implicitly assumed by Peirce (but this specific assumption
is not strictly needed for the main conclusions of the paper), 
we judge that the value of $p$ could lie with equal probability
in any small interval in the range between 0 and 1 
(`uniform prior'). 
Therefore we can use the results obtained in the 
previous section. 

Let us now calculate the weight of evidence provided by the observation
of black or white in favor of the occurrence of black or white. 
We need to calculate the Bayes factor that changes the odd ratio
$P(\mbox{W}\,|\,\mvec{E},B_?,I)/P(\mbox{B}\,|\,\mvec{E},B_?,I)$ 
if we add a further observation $E_{n+1}$, schematically
\begin{eqnarray}
\frac{P(\mbox{W}\,|\,\mvec{E},B_?,I)}{P(\mbox{B}\,|\,\mvec{E},B_?,I)} 
&\rightarrow & \frac{P(\mbox{W}\,|\,\mvec{E},E_{n+1},B_?,I)}
                    {P(\mbox{B}\,|\,\mvec{E},E_{n+1}, B_?,I)}\,.
\end{eqnarray}
This updating factor cannot be calculated directly in an easy way, 
but it can be nevertheless valuated indirectly 
by its definition ('final odds
divided initial odds'). In fact, the evaluation of initial and final 
odds is very simple, just applying Laplace's rule. 
For the former we have
\begin{eqnarray}
O_{W,B}(n_W,n_B,I) &=&\frac{P(\mbox{W}\,|\,\mvec{E},B_?,I)}
                {P(\mbox{B}\,|\,\mvec{E},B_?,I)}  
           = \frac{n_W+1}{n_B+1}\,.  \label{eq:odd0}
\end{eqnarray}
The observation of a new white or of a new black changes this
ratio in 
\begin{eqnarray}
O_{W,B}(n_W+1,n_B,I) &=& \frac{(n_W+1)+1}{n_B+1}\label{eq:odd0+w}\\
O_{W,B}(n_W,n_B+1,I) &=& \frac{n_W+1}{(n_B+1)+1}\label{eq:odd0+b}\,,
\end{eqnarray}
respectively.
Dividing Eqs.~(\ref{eq:odd0+w}) and (\ref{eq:odd0+b}) by 
(\ref{eq:odd0}) we get the 
updating factors of interest:
\begin{eqnarray}
\tilde O_{W,B}(\mbox{W},n_W,n_B,I) &=& \frac{n_W+2}{n_W+1}
\label{eq:BF_W_nW} \\
\tilde O_{W,B}(\mbox{B},n_W,n_B,I) &=& \frac{n_B+1}{n_B+2}\,.
\label{eq:BF_B_nB}
\end{eqnarray}
Contrary to other usual Bayes factors, they depend on the 
{\it previous amount}
 of like-color balls already observed. [We can easily understand
that they also depend on the priors on the box composition,
and therefore Eqs.~(\ref{eq:BF_W_nW})-(\ref{eq:BF_B_nB}) are
only valid for a uniform prior.]

Let us apply these formulae to Peirce's example. The first observation 
of black yields an updating factor of $1/2$, 
or $\Delta$JL of $-0.30$; the second $2/3$, or $\Delta\mbox{JL}=-0.18$;
the third $3/4$, or  $\Delta\mbox{JL}=-0.12$; and so on. 
The updating factor produced by the 20th observation 
is only 20/21=0.952, which corresponds to the little 
weight of evidence $\Delta\mbox{JL}=-0.02$.
The overall factor is 1/21=0.048 
($\Delta\mbox{JL}=-1.32$), which is also equal
to the final odds (the two hypotheses were considered
initially equally likely), from which a probability of 
4.5\% for white and 95.5\% for black can be calculated. 

If the 21-st extraction results, instead, in white, 
the new updating factor is 2 ($\Delta\mbox{JL}=0.30$), 
changing sizeable the overall updating factor, that becomes
then $2/21$, or  $\Delta\mbox{JL}=-1.02$, thus almost
doubling the probability of white, that becomes then
$8.7\%$. This is very interesting: the first time either color
occurs it changes the odds by a factor of two in favor of that color
($|\Delta\mbox{JL}=0.30|$). 
A second observation of white gives an updating factor
of $4/3$ ($\Delta\mbox{JL}=0.12$) and so on. 

If we observe 1010 black and 990 white balls, the updating factor
can be divided into the product of a factor given to white beans and
a factor given to black ones, i.e. 
$\tilde O_{1,2}(990n_W,1010n_B,I) = \tilde O_{1,2}(990n_W,I)\times 
\tilde O_{1,2}(1010n_B,I)$, 
with 
\begin{eqnarray}
 \tilde O_{1,2}(1010n_B,I) &=& \frac{1}{2}\times\frac{2}{3}\times
\frac{3}{4}\times\cdots\times 
\frac{1010}{1011} \\
 \tilde O_{1,2}(990n_W,I) &=& 2\times \frac{3}{2} \times \frac{4}{3}
 \times\cdots\times 
\frac{991}{990}\,.
\end{eqnarray}
In the product the first 990 factors of  
$\tilde O_{1,2}(1010n_B,I)$ are simplified by  
$\tilde O_{1,2}(990n_W,I)$
and the final result is 
\begin{eqnarray}
\tilde O_{1,2}(990n_W,1010n_B,I) &=& \frac{991}{992}\times 
\frac{992}{993} \times \cdots \times \frac{1010}{1011} 
 =  \frac{991}{1011}\,:
\end{eqnarray}
the 20 residual `arguments' in favor of black are not the early
ones (in which case the result would be equal to the observation
of twenty black in a row) but the late ones, individually very small
($\Delta$JL about $-0.0004$ each). All together they provide
a negligible weight of evidence in favor of black, 
with a combined $\Delta$JL of  $-0.0087$, and the 
result of Eq.~(\ref{eq:PB_1010_990}) is reobtained.

\section{Conclusions}
This example shows the danger of drawing quantitative conclusions
from qualitative, intuitive considerations (an issue
extensively discussed in Ref.~\cite{Colombo}). 
Yes, each observation of black is `an argument' in favor
of the opinion that the bean `under the thimble' is black. 
But the arguments do have the same strength and then 
the final `intensity of belief' (to use a very interesting 
expression by Peirce)
 does not depend simply 
on the difference of their numbers in favor of either color.

\end{document}